\documentclass[11pt,reqno]{amsart}
\usepackage{xcolor}
\usepackage{iftex}
\ifPDFTeX
\usepackage[utf8]{inputenc}
\fi
\usepackage[T1]{fontenc}
\usepackage{amsmath,amssymb,amsthm,mathtools}
\usepackage[shortlabels]{enumitem}
\usepackage[margin=1.05in]{geometry}
\usepackage{tikz}
\usepackage[hidelinks]{hyperref}
\usepackage[normalem]{ulem}

\allowdisplaybreaks
\numberwithin{equation}{section}

\newtheorem{theorem}{Theorem}[section]
\newtheorem{proposition}[theorem]{Proposition}
\newtheorem{lemma}[theorem]{Lemma}
\newtheorem{corollary}[theorem]{Corollary}
\theoremstyle{definition}

\theoremstyle{remark}

\newtheorem{example}[theorem]{Example}

\newcommand{\R}{\mathbb R}
\newcommand{\one}{\mathbf 1}

\newcommand{\cL}{\mathcal L}

\title[Prescribed curvature flow on finite trees]
{Prescribed inverse-weight Lin-Lu-Yau curvature flow on finite trees}
\author{Shuliang Bai}
\address{Beijing Institute of Mathematical Sciences and Applications,
Beijing, China}
\email{baishuliang@bimsa.cn}
\author{Haoxuan Cheng}
\address{School of Mathematical Sciences, Fudan University,
Shanghai 200433, China}
\email{hxcheng25@m.fudan.edu.cn}
\author{Bobo Hua}
\address{School of Mathematical Sciences, LMNS, Fudan University,
Shanghai 200433, China}
\email{bobohua@fudan.edu.cn}
\date{}

\begin{document}

\begin{abstract}

We characterize the image of the inverse-weight Lin–Lu–Yau curvature map on finite (edge-weighted) trees. We prove that the prescribed curvature flow converges if and only if the target curvature is attainable; in that case, the convergence rate is exponential. The key ingredient is to formulate the prescribed curvature flow as the gradient flow of a convex functional.


\end{abstract}
\maketitle

\section{Introduction}
Assign positive lengths (weights) to the edges of a finite tree, and let these lengths
determine both the path metric and a random walk whose transition probabilities
are proportional to the inverse edge lengths.  We consider two questions:
which vectors can arise as Lin--Lu--Yau edge curvatures of such metrics,
and when does the associated prescribed-curvature flow converge to a positive metric?

Ollivier \cite{Ollivier2009} introduced Ricci curvature through optimal
transport.  Lin, Lu, and Yau \cite{LinLuYau2011} defined the corresponding
limiting curvature as the idleness tends to one.
Jost and Liu \cite{JostLiu2014} related curvature to local graph structure;
M\"unch and Wojciechowski \cite{MunchWojciechowski2019} developed a
Laplacian formulation for weighted graphs.  Dependence on idleness was
studied by Bourne, Cushing, Liu, M\"unch, and Peyerimhoff
\cite{BourneEtAl2018} for combinatorial graphs, and the large-idleness
relation between Ollivier and Lin--Lu--Yau curvature was extended to
weighted graphs by Liu and Yang \cite{LiuYang2026}.

For continuous-time graph Ricci flow, Bai, Lin, Lu, Wang, and Yau
\cite{BaiEtAl2024} established existence and uniqueness results and
analyzed paths and stars.  Li and M\"unch \cite{LiMunch2026} treated a
discrete-time flow through nonlinear Markov chains, while Ma and Yang
\cite{MaYang2025} studied piecewise-linear flows with surgery.

The inverse-weight Ricci flow on finite trees was studied by Bai, Hua, Lin,
and Liu \cite{BaiHuaLinLiu2025}.  More recently, Bai and Hua \cite{BaiHua2026}
proved that the curvature along the unprescribed flow converges to the
minimum-norm point of a polyhedron determined by the tree.  They also obtained
a convex potential for the curvature map and described the freedom to rescale
positive constant-curvature solutions on conductive components.  In particular,
\cite[Theorem~1.6]{BaiHua2026} gives a necessary and sufficient condition,
in terms of strict subset inequalities, for the existence of a positive
constant-curvature solution on a connected conductive graph.

In this paper, we consider arbitrary prescribed edge-curvature vectors.
We characterize the full image of the curvature map on positive metrics by
component equalities and strict subset inequalities.  We then show that the
prescribed-curvature flow converges to a positive metric precisely when the
target belongs to this image.  For such targets, convergence is exponential,
and the limiting metric is uniquely determined by the initial products of
edge lengths on the conductive components.

\pagebreak[3]
Let $T=(V,E)$ be a finite simple, undirected tree with vertex set $V$,
edge set $E$, and at least one edge.  Write $\R^E$ for the space of real
vectors indexed by $E$, and $\R_+^E$ for its subset of vectors with strictly
positive coordinates.  Assign a positive length $w_e$ to each edge $e$ and
write $w=(w_e)_{e\in E}\in\R_+^E$ for the edge-length vector.
These lengths determine the path metric $d_w$ on $V$.  At a vertex
$x$, let $d_x$ denote its combinatorial degree and put
\[
 D_x(w):=\sum_{e\ni x}\frac1{w_e}.
\]
For $\alpha\in[0,1)$, the lazy probability measure based at $x$ is given by
\[
 m_x^\alpha(x)=\alpha,
 \qquad
 m_x^\alpha(z)=(1-\alpha)\frac{w_{xz}^{-1}}{D_x(w)}
 \quad(z\sim x).
\]
The measure vanishes elsewhere.  Shorter edges carry larger transition
probabilities: in electrical-network terms, $w_e$ is the resistance and
$w_e^{-1}$ the conductance.  For $e=xy$, the
Lin--Lu--Yau curvature is
\[
 \kappa_e(w)
 =\lim_{\alpha\uparrow1}
 \frac{1-W_1(m_x^\alpha,m_y^\alpha)/w_e}{1-\alpha},
\]
where $W_1$ is the Wasserstein distance associated with $d_w$, namely
the minimum transportation cost between the two measures for cost $d_w$;
see \cite[Chapter~1]{Villani2003}.  Varying a length changes both the
transportation cost and the transition kernel, making the curvature
nonlinear even on a tree.

Write $\kappa(w)=(\kappa_e(w))_{e\in E}\in\R^E$ for the corresponding
edge-curvature vector.  The curvature map sends the edge-length vector
to this curvature vector:
$$\kappa:\R_+^{E}\to \R^{E}, w\mapsto \kappa(w).$$

The curvature inequalities obtained in \cite{BaiHua2026} define a polyhedron
containing all curvature vectors of positive metrics.  Our first result
identifies the part of this polyhedron that is actually attained.

For a vertex $x$, let $E(x)$ be the set of edges incident to $x$.  For
$A\subseteq E$, set
\begin{equation}\label{eq:betalambda}
 \beta(A):=\sum_{\substack{x\in V\\E(x)\subseteq A}}(2-d_x).
\end{equation}
This quantity sums the vertex contributions $2-d_x$ over vertices whose
incident edges all lie in $A$, and depends only on the combinatorial structure
of the tree and the edge set $A$.  Proposition~\ref{prop:subset} shows that
$\beta(A)$ is an upper bound for the total curvature on $A$.
For $z\in\R^E$, write $z(A)=\sum_{e\in A}z_e$, and define
\begin{equation}\label{eq:basepolyhedron}
 \mathcal{K}:=
 \left\{z\in\R^E:z(A)\le\beta(A)\text{ for every }A\subseteq E,
                         \ z(E)=\beta(E)\right\}.
\end{equation}

{
Let $G_E$ be the graph whose vertex set is $E$, where two distinct edges are
adjacent if they meet at a vertex of degree at least three.  We call $G_E$ the
\emph{conductive line graph}.  Let
\[
                         E=C_1\sqcup\cdots\sqcup C_m
\]
be the decomposition into the connected components of $G_E$.  Define
\begin{equation}\label{eq:strict-image}
 \mathcal{K}':=
 \left\{z\in\R^E:
 \begin{array}{ll}
 z(C_i)=\beta(C_i),&1\le i\le m,\\
 z(A)<\beta(A),&\emptyset\ne A\subsetneq C_i,\quad 1\le i\le m
 \end{array}\right\}.
\end{equation}
}

\begin{theorem}\label{thm:image-general}

The image of the curvature map on positive edge weights is
\[
          \kappa(\R_+^E)=\mathcal{K}'.
\]

\end{theorem}

For a target $\kappa^*\in\R^E$, we consider the prescribed curvature flow
\begin{align}\label{eq:normalized}
 \dot w_e=-w_e\bigl(\kappa_e(w)-\kappa_e^*\bigr).
\end{align}
A dot denotes differentiation with respect to time $t$.

{
\begin{theorem}\label{thm:flow-convergence}
For every $\kappa^*\in\R^E$ and $w(0)\in\R_+^E$, the prescribed curvature
flow \eqref{eq:normalized} with initial value $w(0)$ has a unique positive
solution for all $t\ge0$.
This solution converges to a positive weight $w^\infty$ if and only if
$\kappa^*\in\mathcal K'$.  In this case $w^\infty$ is uniquely determined by
\[
 \kappa(w^\infty)=\kappa^*,\qquad
 \prod_{e\in C_i}w_e^\infty=\prod_{e\in C_i}w_e(0),
 \quad 1\le i\le m.
\]
There are constants $C,c>0$, depending on $T$, $\kappa^*$, and $w(0)$,
such that, for all $t\ge0$,
\[
 \|\log w(t)-\log w^\infty\|+\|w(t)-w^\infty\|
       +\|\kappa(w(t))-\kappa^*\|\le Ce^{-ct}.
\]
\end{theorem}
}

{
When all internal vertices have degree at least three, $G_E$ is connected,
and the target conditions reduce to
\[
 z(E)=2,\qquad z(A)<\beta(A)\quad(\emptyset\ne A\subsetneq E).
\]
In general, degree-two vertices separate independent groups of variables,
each with its own scale.  Uniqueness therefore requires fixing the scale
on each conductive component.

The prescribed flow considered here differs from flows with a fixed transport
metric or a fixed transition kernel.  For example, Lin and Liu \cite{LinLiu2026}
studied prescribed curvature with fixed graph distances and transition
probabilities proportional to the evolving weights, while Bai, Liu, and Lai
\cite{BaiLiuLai2026} considered an evolving metric with fixed transition
probabilities.  Li, Wang, Yuan, and Zheng \cite{LiWangYuanZheng2026} studied
prescribed Calabi flows for the fixed-distance curvature.  With uniform
neighbor probabilities, Bai, Cheng, and Hua \cite{BaiChengHua2026}
characterized Einstein metrics by a Ricci matrix.
In our setting, the path metric and the inverse-weight transition kernel
both depend on the edge lengths.

Our argument uses the convex potential and component scales of
\cite{BaiHua2026}, following the normalization and convexity
approach of Chow and Luo \cite{ChowLuo2003}.
Related variational formulations include combinatorial Yamabe
flow \cite{Luo2004} and circle-pattern functionals
\cite{BobenkoSpringborn2004}.  For background on base polyhedra and
gradient images of logarithms of finite sums of exponentials, see
\cite{Fujishige1980,Fujishige2005} and
\cite[Theorem~3.3]{WainwrightJordan2008}, respectively.
For the image theorem, fixing the component scales makes the
Hessian positive definite on the normalized space and hence
gives injectivity.  To prove properness, we show that a sequence
escaping to infinity in normalized logarithmic coordinates
has a subsequence along which a proper-subset constraint
becomes an equality in the limit.  Invariance of domain then
gives surjectivity onto $\mathcal K'$.
The image theorem supplies the equilibrium needed for
Theorem~\ref{thm:flow-convergence}.

For an attainable target, the conserved component products fix an affine
space in logarithmic coordinates.  On this space the flow is the negative
gradient flow of the potential, and strict convexity gives a unique
equilibrium.  The image theorem gives coercivity, so the trajectory stays
in a compact sublevel set where the restricted Hessian has a uniform
positive lower bound.  The
gradient-flow estimate yields exponential convergence, first for the
logarithmic coordinates and then for the lengths and curvatures.
Example~\ref{ex:star} exhibits a boundary target on a three-leaf star
for which the curvature converges while the metric degenerates.

Section~\ref{sec:setup} establishes the curvature constraints, and
Section~\ref{sec:heat} gives the convex potential.
Sections~\ref{sec:image} and~\ref{sec:flow} prove the two main results.
Section~\ref{sec:examples} illustrates the role of strict inequalities and
degree-two vertices.
}

\section{Curvature formula and subset inequalities}\label{sec:setup}
The curvature formula gives the inequalities satisfied by attainable targets.
The equality cases determine which components can be scaled independently.

We use the notation of the introduction.  Write $\one$ for the all-one vector
in $\R^E$; all inner products and norms are Euclidean.
The following tree formula is obtained by setting $\gamma(s)=s^{-1}$ in
\cite[Section~2.2]{BaiHuaLinLiu2025}; see also
\cite[Proposition~2.1]{BaiHua2026}.

\begin{proposition}\label{prop:treeformula}
For every edge $e=xy$,
\begin{equation}\label{eq:kappa}
  \kappa_e
  =\frac{2-d_x}{w_eD_x}+\frac{2-d_y}{w_eD_y}.
\end{equation}
In particular, curvature is invariant under the common rescaling $w\mapsto sw$.
Moreover,
\[
 |\kappa_{xy}|\le |2-d_x|+|2-d_y|,
\]
so the curvatures are bounded by a constant depending only on the tree.
\end{proposition}

\begin{proposition}\label{prop:GB}
For every positive metric on a finite tree,
\begin{equation}\label{eq:GB}
             \sum_{e\in E}\kappa_e=2.
\end{equation}
\end{proposition}

{
\begin{proof}[Proof of Proposition~\ref{prop:GB}]
Summing \eqref{eq:kappa} by vertices gives
\[
 \sum_{e\in E}\kappa_e
 =\sum_{x\in V}(2-d_x)
       \frac{\sum_{e\ni x}w_e^{-1}}{D_x(w)}
 =2|V|-2|E|=2.
\]
The last equality uses that $T$ is a tree.  In particular, $\beta(E)=2$.
\end{proof}
}
The following subset bound and equality criterion follow from
\cite[proof of Theorem~1.3]{BaiHua2026}.  Equality identifies the independent components, and the strict inequalities characterize
the curvatures of positive metrics.
\begin{proposition}\label{prop:subset}
For every nonempty subset $A\subseteq E$,
\[
 \sum_{e\in A}\kappa_e\le\beta(A).
\]
The inequality is strict if and only if there is a vertex $x\in V$ such that
$d_x\ge3$, $E(x)\cap A\ne\emptyset$, and $E(x)\not\subseteq A$; otherwise
equality holds.
\end{proposition}
\begin{proof}
    Let $A\subseteq E$ be nonempty. Summing the curvature formula over $A$ and
grouping the terms by vertices gives
\begin{equation}\label{eq:subsetcurvature}
 \sum_{e\in A}\kappa_e
 =\sum_{\substack{x\in V\\A\cap E(x)\ne\varnothing}}
 (2-d_x)
 \frac{\displaystyle\sum_{e\in A\cap E(x)}w_e^{-1}}
      {\displaystyle\sum_{f\in E(x)}w_f^{-1}}.
\end{equation}

 If $E(x)\subseteq A$, its contribution equals $2-d_x$ and
is already included in $\beta(A)$. A vertex meeting both
$A$ and $E\setminus A$ has degree at least two, since a
leaf has only one incident edge. Such a vertex contributes
zero if $d_x=2$.

Consequently, subtracting \eqref{eq:subsetcurvature}
from the definition of $\beta(A)$ gives
\begin{equation}\label{eq:subset-gap}
\begin{aligned}
\beta(A)-\sum_{e\in A}\kappa_e(w)
={}&
\sum_{\substack{x\in V\\d_x\ge3\\
0<|A\cap E(x)|<d_x}}
(d_x-2)\\
&\quad\cdot
\frac{\sum_{e\in A\cap E(x)}w_e^{-1}}{D_x(w)}.
\end{aligned}
\end{equation}

Every summand is strictly positive. Therefore
\begin{equation}\label{eq:subsetbound}
\sum_{e\in A}\kappa_e\le\beta(A),
\end{equation}
with strict inequality precisely when some vertex of
degree at least three is incident to edges on both
sides of the cut $(A,E\setminus A)$.

\end{proof}

\begin{lemma}\label{lem:component-polyhedron}
For every $A\subseteq E$,
\begin{equation}\label{eq:beta-components}
\beta(A)=\sum_{i=1}^m\beta(A\cap C_i).
\end{equation}
Moreover, with
\[
H:=\bigl\{z\in\R^E:
z(C_i)=\beta(C_i),\ 1\le i\le m\bigr\},
\]
one has
\[
\operatorname{aff}(\mathcal K)=H,
\qquad
\operatorname{ri}(\mathcal K)=\mathcal K'.
\]
Here $ri$  denotes relative interior in the affine hull.
\end{lemma}
\begin{proof}
 If $d_x\ge3$, all edges in $E(x)$ belong to one conductive component.
 The same holds for a leaf, while $2-d_x=0$ at a degree-two vertex.
 Splitting the vertex sum defining $\beta(A)$ therefore gives
 \eqref{eq:beta-components}.
 
 For $z\in\mathcal K$, the inequalities $z(C_i)\le\beta(C_i)$,
 together with $z(E)=\beta(E)=\sum_i\beta(C_i)$, force
 $z(C_i)=\beta(C_i)$ for every $i$. Thus $\mathcal K\subseteq H$.
 Conversely, if $z\in H$ satisfies
 $z(B)<\beta(B)$ for every nonempty proper subset $B\subsetneq C_i$,
 then, for any $A\subseteq E$, applying these inequalities to
 $A\cap C_i$ (with equality when the intersection is empty or all of
 $C_i$) gives
 \[
  z(A)=\sum_i z(A\cap C_i)
  \le\sum_i\beta(A\cap C_i)=\beta(A).
 \]
 In particular, $\mathcal K'$ is a relatively open subset of $H$
 contained in $\mathcal K$.
 
 The set $\mathcal K'$ is nonempty: for any positive metric $w$,
 Proposition~\ref{prop:subset} gives the component equalities, and the
 connectedness of each $C_i$ in $G_E$ makes every proper nonempty
 subset inequality strict.  Thus $\operatorname{aff}(\mathcal K)=H$
 and $\mathcal K'\subseteq\operatorname{ri}(\mathcal K)$.
 Conversely, each local inequality $z(B)\le\beta(B)$ with
 $\emptyset\ne B\subsetneq C_i$ is strict at $\kappa(w)$, so it is
 not an equality throughout $\mathcal K$. Such an inequality must be
 strict at every point of the relative interior. Hence
 $\operatorname{ri}(\mathcal K)\subseteq\mathcal K'$.
 \end{proof}

\section{Convex potential}\label{sec:heat}
In logarithmic coordinates, the prescribed flow is a Euclidean gradient flow.
Adding a linear target term to the potential in
\cite[Section~3.2]{BaiHua2026} leaves its Hessian unchanged.
The kernel corresponds to independent rescaling of the components.
Let $\kappa^*\in\R^E$.
Write $y_e=\log w_e$ and define
\begin{equation}\label{eq:Phi}
 \Phi(y):=\sum_{x\in V}(d_x-2)\log D_x(y)-\sum_{e\in E}\kappa_e^* y_e,
 \qquad
 D_x(y):=\sum_{e\ni x}e^{-y_e}.
\end{equation}

{
For $z\in\R^E$, define the linear operator $\cL(w):\R^E\to\R^E$ by
\begin{equation}\label{eq:laplacian}
 (\cL(w)z)_e
 =\sum_{\substack{x\in e\\d_x\ge3}}
    \sum_{\substack{f\ni x\\f\ne e}}
       \frac{d_x-2}{D_x(w)^2w_ew_f}(z_e-z_f).
\end{equation}
Here $x\in e$ means that $x$ is an endpoint of $e$.
This is the weighted graph Laplacian on $G_E$.
}
\begin{proposition}\label{prop:potential}
For every $e\in E$,
\begin{equation}\label{eq:gradPhi}
                    \frac{\partial\Phi}{\partial y_e}=\kappa_e-\kappa_e^*,
 \qquad             \nabla^2\Phi(y)=\cL(w)\ge0.
\end{equation}
The function $\Phi$ is convex and the prescribed curvature flow is
$\dot y=-\nabla\Phi,$ the negative Euclidean gradient flow in $y$.  
The kernel of $\nabla^2\Phi$ consists of the vectors $v$ such that $v_e=v_f$
whenever $e$ and $f$ meet at a vertex of degree at least three.  Equivalently,
$v$ is constant on each connected component of the conductive line graph.  In particular,
if the conductive line graph is connected, then $\Phi$ is strictly convex on
$\one^\perp$.
\end{proposition}

{
\begin{proof}
Since $D_x(y)=\sum_{f\ni x}e^{-y_f}$, differentiation gives
\[
 \frac{\partial\Phi}{\partial y_e}
 =\sum_{x\in e}(d_x-2)\frac{-e^{-y_e}}{D_x(y)}-\kappa_e^*
 =\sum_{x\in e}\frac{2-d_x}{w_eD_x(w)}-\kappa_e^*
 =\kappa_e-\kappa_e^*.
\]
For $v\in\R^E$, a second differentiation gives
\[
 v^T\nabla^2\Phi(y)v
 =\sum_{x:\,d_x\ge3}(d_x-2)
 \left[
  \frac{\sum_{e\ni x}e^{-y_e}v_e^2}{D_x(y)}
  -\left(\frac{\sum_{e\ni x}e^{-y_e}v_e}{D_x(y)}\right)^2
 \right]\ge0.
\]
The off-diagonal Hessian entry for distinct edges meeting at $x$ is
$-(d_x-2)/(D_x(w)^2w_ew_f)$ when $d_x\ge3$, and is zero otherwise;
each diagonal entry is the negative sum of the off-diagonal entries in its row.
Thus the Hessian is $\cL(w)$.
Each expression in brackets above is a weighted variance.  Since all weights
are positive, it vanishes if and only if $v$ is constant on $E(x)$.  The leaf
terms are linear and the degree-two terms vanish.  Hence the whole quadratic
form vanishes if and only if $v$ is constant on every connected component of
the conductive line graph.  Finally,
$\dot y_e=\dot w_e/w_e=-(\kappa_e-\kappa_e^*)$, so the prescribed curvature
flow is $\dot y=-\nabla\Phi$.
\end{proof}
}

\section{The image of the curvature map}\label{sec:image}
Fixing the component scales makes the curvature map injective.  To obtain
Theorem~\ref{thm:image-general}, we prove that the normalized map is proper
and apply the following consequence of invariance of domain.

A map $f:\R^n\to V$ with open subset $V\subset \R^n$ is called \emph{proper} if $f^{-1}(K)$ is compact for any compact $K\subset V.$
\par\noindent\begin{minipage}{\linewidth}
\begin{corollary}\label{cor:proper-injective}
   Let \(V\subset\mathbb R^n\) be a connected open set. If

$$ f:\mathbb R^n\to V $$

is continuous, injective, and proper, then \(f\) is a homeomorphism from \(\mathbb R^n\) onto \(V\).
\end{corollary}
\end{minipage}\par
\begin{proof}
    By the invariance of domain theorem, $f(\R^n)$ is open in $V.$ The image of a proper map into a locally compact Hausdorff space is closed.  Hence $f(\R^n)=V.$
\end{proof}

We first prove the image theorem when no internal vertex has degree two.
The same argument will then be applied on each conductive component.
\begin{theorem}\label{thm:image-no-degree-two} If all internal vertices have degree at least 3, then
    $\kappa(\R_+^E)=\mathcal{K}'.$
\end{theorem}
\begin{proof}
  Consider the curvature map on $\R_+^E.$
The subset inequalities give $$\kappa:\R_+^E\to \mathcal{K}'.$$ %
With $u=\log w,$ or equivalently $w=e^u,$ we define
$$\hat{\kappa}:\Omega\to \mathcal{K}', \hat{\kappa}(u)_e=\kappa(e^u)_e,$$
where $\Omega=\{u\in \R^E: \sum_{e\in E}u_e=0\}\simeq \R^{|E|-1}.$ Note that $\dim \Omega=\dim \mathcal{K}'=|E|-1.$

By Proposition~\ref{prop:potential}, the restriction of $\Phi$ to $\Omega$
is strictly convex, so $\hat{\kappa}$ is injective.  Here $u$ denotes the
same logarithmic coordinate as $y$ in Section~\ref{sec:heat}, restricted
to the normalized space $\Omega$.

To prove that $\hat{\kappa}$ is proper, let $K\subset \mathcal{K}'$ be compact.
We must show that $\hat\kappa^{-1}(K)$ is compact.  This preimage is closed, so if it were not compact, there would
be a sequence $\{u_i\}\subset\hat\kappa^{-1}(K)$ such that $\lim_i|u_i|=\infty.$ Passing to a subsequence, every coordinate $u_i(e)$ has a limit in the extended real line. Since $\sum_eu_i(e)=0$ and $|u_i|\to\infty$, at least one coordinate tends to $-\infty$ and at least one does not.  Set
\[
 A:=\left\{e\in E:\lim_{i\to\infty}u_i(e)>-\infty\right\}.
\]
Then $\emptyset\ne A\subsetneq E$.  For $w_i=e^{u_i}$, the sequence
$w_i(e)$ is bounded away from zero for $e\in A$, while $w_i(e)\to0$ for
$e\notin A$.

At a vertex $x$ meeting both $A$ and $E\setminus A$, the numerator in the $A$-ratio in \eqref{eq:subsetcurvature} stays bounded, while the denominator tends to infinity.  If $E(x)\subseteq A$, that ratio is identically one.  Hence, by \eqref{eq:subsetcurvature},
$$\lim_{i\to \infty}\sum_{e\in A}\kappa_e(w_i)=\beta(A).$$
But $\hat{\kappa}(\{u_i\})={\kappa}(\{w_i\})\subset K$, where $K$ is a compact subset of $\mathcal{K}'.$
Compactness of $K$ gives a positive lower bound for $\beta(A)-z(A)$ on $K$,
whereas the preceding limit makes $\beta(A)-\hat\kappa(u_i)(A)$ tend to zero.
This contradiction proves that $\hat\kappa$ is proper.
After identifying $\Omega$ and the affine hyperplane $\{z:z(E)=2\}$ with $\R^{|E|-1}$,
Corollary~\ref{cor:proper-injective} implies that $\hat\kappa$ maps $\Omega$ onto $\mathcal{K}'.$
 
\end{proof}

{
We now prove Theorem~\ref{thm:image-general} for an arbitrary finite tree.
A vertex of degree two separates its two incident edges into
different components of $G_E$, whereas all edges incident to a vertex of
degree at least three belong to the same component.

\begin{proof}[Proof of Theorem~\ref{thm:image-general}]
For the general tree, set
\[
 \Omega:=\left\{u\in\R^E:
       \sum_{e\in C_i}u_e=0,\ 1\le i\le m\right\}
\]
which agrees with the preceding definition of $\Omega$ when $G_E$ is
connected.  Both $\Omega$ and the affine subspace
\[
 \left\{z\in\R^E:z(C_i)=\beta(C_i),\ 1\le i\le m\right\}
\]
have dimension $|E|-m$, and $\mathcal{K}'$ is a convex open subset of this
affine subspace.

We first show that $\kappa(\R_+^E)\subseteq\mathcal{K}'$.  For each component
$C_i$, equality holds in \eqref{eq:subsetbound}, since every vertex meeting
both $C_i$ and $E\setminus C_i$ has degree two.  Thus
$\kappa(w)(C_i)=\beta(C_i)$.  If $\emptyset\ne A\subsetneq C_i$, the connectedness
of $C_i$ in $G_E$ gives a vertex of degree at least three which is incident to
an edge in $A$ and an edge in $C_i\setminus A$.  The inequality in
\eqref{eq:subsetbound} is therefore strict.

Consider the restriction
\[
             \hat{\kappa}:\Omega\longrightarrow\mathcal{K}',
             \qquad \hat{\kappa}(u)=\kappa(e^u).
\]
By Proposition~\ref{prop:potential}, the kernel of the Hessian of $\Phi$ is
the set of vectors which are constant on each $C_i$.  If such a vector $v$
belongs to $\Omega$, then, for any $e\in C_i$,
$0=\sum_{f\in C_i}v_f=|C_i|v_e$.  Thus $v=0$.  Hence the
restriction of $\Phi$ to $\Omega$ is strictly convex.  Since
$\nabla\Phi(u)=\hat{\kappa}(u)-\kappa^*$, the map $\hat{\kappa}$ is injective.

We next prove that $\hat{\kappa}$ is proper.  Let $K\subset\mathcal{K}'$ be compact.
The set $\hat{\kappa}^{-1}(K)$ is closed in $\Omega$.  Suppose that it is not
compact.  Then there is an unbounded sequence
$\{u_i\}\subset\hat{\kappa}^{-1}(K)$.  For some $C_j$, after passing to a subsequence,
$\|u_i|_{C_j}\|\to\infty$.  Passing to a further subsequence, every coordinate
has a limit in the extended real line on $C_j$.  Since
$\sum_{e\in C_j}u_i(e)=0$, the set
\[
       A:=\Bigl\{e\in C_j:\lim_{i\to\infty}u_i(e)>-\infty\Bigr\}
\]
is a nonempty proper subset of $C_j$: at least one coordinate tends to
$-\infty$, while at least one does not.  If a vertex $x$ of degree at least
three meets both $A$ and $E\setminus A$, then all edges in $E(x)$ belong to
$C_j$.  The numerator of the $A$-ratio in
\eqref{eq:subsetcurvature} stays bounded, while its denominator tends to
infinity.  A vertex of degree two contributes zero, and if $E(x)\subseteq A$,
the ratio is identically one.  Therefore
\[\lim_{i\to\infty}\hat{\kappa}(u_i)(A)=\beta(A).
\]
On the other hand, compactness of $K$ and the strict inequality defining
$\mathcal{K}'$ give
\[
       \min_{z\in K}\bigl(\beta(A)-z(A)\bigr)>0,
\]
a contradiction.  Thus $\hat{\kappa}$ is proper.  If $|E|-m=0$, then $\Omega$ and $\mathcal{K}'$ are both singletons, and the conclusion follows.  Otherwise, after identifying $\Omega$ and the preceding affine subspace with $\R^{|E|-m}$,
Corollary~\ref{cor:proper-injective} implies that
$\hat{\kappa}(\Omega)=\mathcal{K}'$.  Since the curvature is unchanged when all weights in any one component $C_i$ are multiplied by the same positive constant, every positive weight vector can be normalized into $\Omega$ by multiplying the weights in $C_i$ by
$\exp\bigl(-|C_i|^{-1}\sum_{e\in C_i}\log w_e\bigr)$.
Indeed, every vertex term with nonzero coefficient involves edges from only one component, while a vertex of degree two contributes zero.  Hence
$\kappa(\R_+^E)=\mathcal{K}'$.
\end{proof}

}

\section{Convergence of the prescribed flow}\label{sec:flow}
Bounded curvature gives global existence for every target.  For an
attainable target, the image theorem supplies an equilibrium with the
prescribed component products.  Coercivity and the Hessian bound then give
exponential convergence to this equilibrium.

{
\begin{proof}[Proof of Theorem~\ref{thm:flow-convergence}]
Put $y_e=\log w_e$.  By Proposition~\ref{prop:potential},
\[
                  \dot y=-\nabla\Phi(y)=-(\kappa(w)-\kappa^*).
\]
The curvature is uniformly bounded by Proposition~\ref{prop:treeformula}.
Thus the right-hand side is smooth and bounded on $\R^E$, and the standard
continuation theorem gives a unique solution for all $t\ge0$.  Hence
$w_e(t)=e^{y_e(t)}>0$ at every finite time.

Suppose that $\kappa^*\in\mathcal K'$.  Since
$\kappa(w)(C_i)=\kappa^*(C_i)=\beta(C_i)$, we have
\[
 \frac{d}{dt}\sum_{e\in C_i}y_e(t)=0,
 \qquad 1\le i\le m.
\]
Thus the products $\prod_{e\in C_i}w_e(t)$ are constant, and the initial value
fixes the invariant affine space on each component.  Moreover, the equations
split over the components $C_i$: a vertex of degree two contributes zero,
while all edges incident to a vertex of degree at least three belong to the
same component.
Accordingly, $\Phi$ is a sum of functions in the disjoint sets of variables
$y|_{C_i}$.

Fix $C_i$ and restrict the corresponding component of $\Phi$ to the affine
space determined by the initial value,
\[
 \left\{y|_{C_i}:\sum_{e\in C_i}y_e=
                         \sum_{e\in C_i}y_e(0)\right\}.
\]
If $|C_i|=1$, then $\kappa_e(w)=\kappa_e^*=\beta(C_i)$ on its unique edge,
so the flow on this component is constant.  Suppose that $|C_i|>1$.
Theorem~\ref{thm:image-general} gives a positive weight realizing $\kappa^*$.  By a
common rescaling on $C_i$, this weight can be chosen in the displayed affine
space without changing its curvature.  It is therefore a critical point of
the restricted potential.  Proposition~\ref{prop:potential} shows that the
restricted Hessian is positive definite, so the critical point is unique.

We next prove that the restricted potential is coercive.  In what follows,
all gradients and Hessians are taken on the displayed affine space.  Its
tangent space is
\[
             \left\{v\in\R^{C_i}:\sum_{e\in C_i}v_e=0\right\}.
\]
The restricted gradient is
$\kappa(w)|_{C_i}-\kappa^*|_{C_i}$, whose coordinates sum to zero.  Since all
the inequalities defining $\mathcal K'$ are strict at $\kappa^*$, every
sufficiently small tangent vector $q$ keeps
$\kappa^*|_{C_i}+q$ within all the inequalities in
\eqref{eq:strict-image} associated with $C_i$.
Keeping the target fixed on the other components and applying
Theorem~\ref{thm:image-general}, followed by a rescaling on $C_i$, shows that
$q$ belongs to the restricted gradient image on the displayed affine space.
Thus zero is an interior point of this image in the tangent space.  In
particular, there is $r>0$ such that $rv$ belongs to the gradient image for
every unit tangent vector $v$.  Let $p_v$ be the unique point with
$\nabla\Phi(p_v)=rv$.
Translation by a vector which is constant on $C_i$ identifies the displayed
affine space with the zero-sum space used in the proof of
Theorem~\ref{thm:image-general} and does not change curvature.  The
componentwise properness argument there shows
that the set of all such $p_v$ is compact.  For any $y$ in the affine space, put
\[
 h=y-y(0)|_{C_i},\qquad v=\frac{h}{|h|}
\]
when $h\ne0$.  The supporting inequality for the convex function $\Phi$ gives
\begin{align*}
 \Phi(y)
 &\ge \Phi(p_v)+\langle\nabla\Phi(p_v),y-p_v\rangle\\
 &=r|h|+\Phi(p_v)+\langle rv,y(0)|_{C_i}-p_v\rangle.
\end{align*}
The last two terms have a uniform lower bound because the points $p_v$ form a
compact set.  Thus $\Phi(y)\ge r|h|-C$ for some constant $C$, which proves
coercivity.

It remains to apply the standard finite-dimensional convex-gradient-flow
argument, analogous to \cite[Propositions~3.7--3.8 and \S6.2]{ChowLuo2003};
we give the estimates here.  Along
the flow,
\[
             \frac{d}{dt}\Phi(y(t))=-|\nabla\Phi(y(t))|^2\le0.
\]
The trajectory on this component therefore stays in a compact sublevel set.
On this set the restricted Hessian is positive definite and continuous, so
there is a constant $\mu>0$ such that
\[
             v^T\nabla^2\Phi(y)v\ge\mu|v|^2
\]
for every tangent vector $v$.  Let $y^\infty$ denote the unique critical point
on this affine space.  By convexity, the segment from $y^\infty$ to $y(t)$
also lies in the same sublevel set.  Integrating the Hessian along this segment
gives
\[
 \frac{d}{dt}\frac12|y(t)-y^\infty|^2
 =-\langle y(t)-y^\infty,\nabla\Phi(y(t))\rangle
 \le-\mu|y(t)-y^\infty|^2.
\]
Gronwall's inequality yields
\[
 |y(t)-y^\infty|\le e^{-\mu t}|y(0)-y^\infty|.
\]
Applying this argument to the finitely many components proves exponential
convergence of $y(t)$.  Since $w=e^y$ and $\kappa(w)$ depend smoothly on $y$
on the compact sublevel set, $w(t)$ and $\kappa(w(t))$ also converge
exponentially.  The limit is the unique $w^\infty$ with the stated curvature
and component products.

Conversely, suppose that $w(t)$ converges to a positive weight $w^\infty$.
Then $y(t)$ converges, while its derivative has the limit
$-(\kappa(w^\infty)-\kappa^*)$.  A convergent differentiable curve cannot have
a nonzero limiting velocity.  Hence this limit is zero and
$\kappa(w^\infty)=\kappa^*$.  Theorem~\ref{thm:image-general} now gives
$\kappa^*\in\mathcal K'$.
\end{proof}
}

\section{Examples and boundary behavior}\label{sec:examples}
{
The star illustrates metric degeneration at a boundary target.  The path
shows why determining the limit requires the edge-length product on each
conductive component.

\begin{example}\label{ex:star}
Let $T$ be the three-leaf star with edges $e_1,e_2,e_3$, and write
$w_j=w_{e_j}$ and $\kappa_j=\kappa_{e_j}$.  Put
\[
 q_j=\frac{w_j^{-1}}{w_1^{-1}+w_2^{-1}+w_3^{-1}},\qquad 1\le j\le3.
\]
The tree formula gives $\kappa_j=1-q_j$.  Conversely, every positive
triple $q$ with $q_1+q_2+q_3=1$ is realized by $w_j=q_j^{-1}$.
Thus the curvature image of positive metrics is
\[
 \{z\in\R^3:z_1+z_2+z_3=2,\ 0<z_j<1\ (1\le j\le3)\}.
\]
The target $(1,1/2,1/2)$ belongs to $\mathcal K$ but is not attained by a
positive metric.

For this target, start the flow with $w_1(0)=a(0)>0$ and
$w_2(0)=w_3(0)=b(0)>0$.  Symmetry and uniqueness preserve
$w_1(t)=a(t)$ and $w_2(t)=w_3(t)=b(t)$.  Direct substitution gives
\[
 \dot a=\frac{ab}{b+2a},\qquad
 \dot b=-\frac{b^2}{2(b+2a)},\qquad
 a(t)b(t)^2=a(0)b(0)^2.
\]
Here $a$ increases and $b$ decreases.  If $a$ were bounded, the product
identity would give positive finite limits for both variables, while
$\dot a$ would have a positive limit, a contradiction.  Hence
$a\to\infty$ and $b\to0$, so $q_1\to0$ and $q_2=q_3\to1/2$.
The curvature converges to the prescribed vector although the metric has
no positive finite limit.
\end{example}

\begin{example}\label{ex:path}
For a path with at least two edges, every conductive component is a
singleton.  The two end edges have curvature $1$, and every other edge
has curvature $0$, independently of the positive lengths.  For a path
consisting of one edge, its curvature is $2$.
In either case $\mathcal K'$ is a singleton.  At its unique attainable
target, the flow is stationary at every initial metric.  The product on
each singleton component is its edge length, so the component products
specify the stationary metric.  A single product over all edges would
not do so for paths with more than one edge.
\end{example}
}

\par\medskip
\noindent\textbf{Acknowledgements.}
S. Bai is supported by NSFC grant no.~12301434. B. Hua is supported by NSFC grant no.~12371056. The author B.H. is grateful to Max Planck Institute for Mathematics in Bonn for its hospitality and financial support.

\par\medskip
\noindent\textbf{Statement of AI Use.}
  ChatGPT was used as a supporting tool for English-language and LATEX editing. The authors take full responsibility for the
 content.

\bibliographystyle{amsplain}
\bibliography{referencesa}

\end{document}